\input amstex
\documentstyle{amsppt}

\topmatter

\title { Character degrees of some $p$-groups}\endtitle
\author {Avinoam Mann}\endauthor
\endtopmatter

It was shown by I.M.Isaacs [Is2], that any set of powers of a prime $p$ which includes 1, can occur as the set of character degrees of some $p$-group. It becomes then of interest to see to what extent that remains true if we consider a particular class of $p$-groups. Isaacs' construction yields groups of nilpotency class 2. Here we consider the other extreme. Recall that a group of order $p^n \ge p^3$ is said to be {\it of maximal class}, if its nilpotency class $cl(G)$ is $n-1$ (see [LGM] for the well developed theory of these groups). Such a group has a non-abelian factor group of order $p^3$, and therefore it has irreducible characters of degree $p$. It was suggested in the last section of [Sl1], that there are further restrictions on the possible character degrees set of a group of maximal class. The present note verifies that conjecture, indeed under weaker assumptions than maximal class, for some results it suffices to assume that the derived subgroup has index $p^2$. This last assumption has quite a few consequences for the structure of the given group, and a secondary aim of this note is to derive some of them (see Theorem 4 and Propositions 5 and 6. These results are going to be applied in [Ma2]). About character degrees our results are:

\proclaim {Theorem 1} Let $G$ be a $p$-group in which $|G:G'| = p^2$. If $G$ has an irreducible character of degree $ > p$, then it has such a character of degree at most $p^{\frac{p+1}{2}}$.\endproclaim

We remark that 2-groups satisfying $|G:G'| = 4$ are of maximal class, and their character degrees are 1 and 2. For all odd primes, however, there exist groups satisfying $|G:G'| = p^2$ which are not of maximal class, and which have irreducible characters of degrees $ > p$.

M.Slattery [Sl2] has constructed, for all primes $p \ge 5$, groups of maximal class whose character degrees are $1, p, p^{\frac{p+1}{2}}$, showing that the bound of Theorem 1 is best possible. It is also easy to see that there exist 3-groups of that type.

\proclaim {Theorem 2} Let $G$ be a $p$-group of maximal class. If $G$ has irreducible characters of
degree higher than $p^2$, then it has such characters of degree at most $p^{\frac{p+3}{2}}$.\endproclaim

\proclaim {Theorem 3} Let $G$ be a $p$-group satisfying $|G:G'| = p^2$. If $|G| \ge p^6$, and $|\gamma_3(G):\gamma_4(G)| > p$ (equivalently, if $|G:\gamma_4(G)| = p^5$), then $G$ has irreducible characters of degree $p^2$.\endproclaim

For the proofs we first quote some of the theory of groups of maximal class, the proofs of which can be found in [LGM] and in [Hu, III.14].
Let $X$ be a $p$-group of maximal class, of class $c$, say. We write $X_i = \gamma_i(X)$ for the terms of the lower central series, for $ 2 \le i \le
c$, and $X_1 = C_X(X_2/X_4)$. These notations will be applied to each group of maximal class that we will
encounter below, e.g. if that group will be denoted by $H$, we will let $H_i$ denote the corresponding subgroups
of $H$, etc. Returning to $X$, we have $X_1 = C_X(X_i/X_{i+2})$, for $2 \le i < c-1$. We call $X_1$ the {\it major
centralizer} of $X$. It is a regular $p$-group, therefore, if $A,B \triangleleft X_1$, then $[A,B^p] = [A,B]^p$
(see [Hu, III.10] for the theory of regular $p$-groups). If $|X| \ge p^{p+2}$, then $X_i^p = X_{i+p-1}$. This
holds also for order $p^{p+1}$, except for one group, the wreath product of two groups of order $p$, in which
$X_i^p = 1$ for all $i$. If $|X| \le p^p$, then $X_i^p = 1$ for $i \ge 2$. If $X_1 \neq C_X(X_{c-1})$, then $X$ is
termed {\it exceptional}. If $X$ is metabelian, it is not exceptional. In a
non-exceptional group, all maximal subgroups different from $X_1$ are of maximal class. Finally, $X$ is of maximal
class iff it contains an element $x\in X$ such that $|C_X(x)| = p^2$.

\smallskip

Next, about groups with a derived subgroup of index $p^2$. Since all normal subgroups of index $p^2$ contain the commutator subgroup, it follows that $G'$ is the only normal subgroup of that index. The factor group $G/G'$ is elementary abelian, and if $x$ and $y$, say, are elements of $G$ that are independent {\it modulo} $G'$, then they generate $G$, and the factor group $G'/\gamma_3(G)$ is generated by the image of $[x,y]$. Therefore $\gamma_3(G)$ has index $p^3$ in $G$ and it is the only normal subgroup of that index. Similarly, $\gamma_3(G)/\gamma_4(G)$ is generated by the images of $[x,y,x]$ and $[x,y,y]$, therefore $\gamma_4(G)$ has index either $p^4$ or $p^5$, and in the first case it is the only normal subgroup of that index.

Recall that $d(X)$ denotes the minimal number of generators of the group $X$.

\proclaim {Theorem 4} Let $G$ be a $p$-group such that $|G| \ge p^3$ and $|G:G'| = p^2$, and let $H$ be a maximal subgroup of $G$. Then

(a) $G/H'$ is of maximal class, and $d(H) \le p$.

(b) If $d(H) < p - 1$, then $|\Phi(H):H'| \le p$, hence $|H:H'|$ equals either  $p^{d(H)+1}$ or $p^{d(H)}$.

(c) If $d(H) = p$, then $G/H' \cong C_pwrC_p$, and $H/H'$ is elementary abelian of order $p^p$.

(d) If $|\gamma_3(G):\gamma_4(G)| = p$, then $|\gamma_4(G):\gamma_5(G)| = p$, and all maximal subgroups $H$ of $G$, save one, satisfy $H' =
\gamma_3(G)$.

(e) If $|\gamma_3(G):\gamma_4(G)| = p^2$, then $H' \neq \gamma_3(G)$ and $H' \not \le \gamma_4(G)$. If $K$ is
another maximal subgroup, then $H' \not\le K'$, and in particular $H' \neq K'$.

(f) If $G$ contains a 2-generator maximal subgroup $H$, then either $H$ is metacyclic, or $G$ contains at most one maximal subgroup $K \neq H$ such that $|K:K'| \ge p^{p+1}$.

(g) If $\gamma_3(G)\neq 1$, then $G^p \le \gamma_3(G)$.
\endproclaim

{\bf Note.} The first claim in (a) was already pointed out in Exercise 1.4 of [Be].

{\bf Proof.} For (a), we may as well assume that $H' = 1$. Then $G = \langle H, x \rangle$ for some element $x$, and
commutation with $x$ induces on $H$ an endomorphism with image $[G,H] = G'$ and kernel $Z(G)$. Since $|H:G'| = p$,
we have $|Z(G)| = p$. Moreover, the same argument shows that all non-abelian quotients of $G$ have centres of
order $p$, hence $G$ is of maximal class.

If $|H:H'| > p^{p-1}$ , then $K := G/H'$ is a $p$-group of maximal class and order at least
$p^{p+1}$. Since $L := H/H'$ is abelian, it is the major centralizer of $K$. It is known that in such groups $K^p
= \gamma_p(K)$, a subgroup of index $p^p$, and that usually $L^p = K^p$, the only exception to this equality
occurring when $K = C_pwr C_p$, when $L$ is elementary abelian of order $p^p$. In the other cases we have $d(H) =
d(L) = \log_p|L:\Phi(L)| = \log_p|L:L^p| = p-1$.

If $|H:H'| \le p^{p-1}$, then $d(H) \le p -1$, and by the above, this is the only case in which strict inequality
is possible. Moreover, in that case $|K| \le p^p$, and it is known that groups of maximal class of these orders
satisfy $|K^p| \le p$, which proves (b).

When $|\gamma_3(G):\gamma_4(G)| = p$, write $C = C_G(G'/\gamma_4(G))$. This is a maximal subgroup, and $C' =
[C,G'] \le \gamma_4(G)$. Let $H$ be another maximal subgroup. By Lemma 30 of [Ma1], $|G':C'H'| \le p$, i.e.
$|G:C'H'| \le p^3$. But $C'H' \le \gamma_3(G)$, and $|G:\gamma_3(G)| = p^3$. Therefore $C'H' = \gamma_3(G)$. Since
$C' \le \gamma_4(G)$, this is possible only if $H' = \gamma_3(G)$.

Since $d(G) = 2$, we have $G'' \le \gamma_5(G)$ [Hu, III.2.12(b)], and since $|G:\gamma_4(G)| = p^4$, an argument of R.Brandl
implies that $G/G''$ is of maximal class (see Theorem C of [Ma1]), which shows that $|\gamma_4(G):\gamma_5(G)| =
p$.

If $|\gamma_3(G):\gamma_4(G)| = p^2$, and $H$ is a maximal subgroup, then $H' \not\le \gamma_4(G)$, because
$G/H'$ is of maximal class. Let $M$ and $N$ be two subgroups lying properly between $\gamma_3(G)$ and
$\gamma_4(G)$, and let $C = C_G(G'/M)$ and $D = C_G(G'/N)$. Then $M$ and $N$ are normal and $C$ and $D$ are
maximal. If $C = D$, then $C' = [C,G'] \le \gamma_4(G)$, contradiction. Thus the $p+1$ subgroups like $M$ and $N$
determine $p+1$ distinct maximal subgroups, i.e. each maximal subgroup is obtained in this way. Moreover, $C' \le
M$, thus $C' \neq \gamma_3(G)$. But $C'D' = \gamma_3(G)$, therefore $D' \not\le C'$.

To prove (f) we suppose that $G$ contains the 2-generator maximal subgroup $H$, and two other maximal subgroups, $K$, and $L$, whose
commutator factor groups have orders at least $p^{p+1}$. If $|H:H^p| \le p^2$, then $H$ is either cyclic or metacyclic. Thus we may assume that $|H:H^p| \ge p^3$. We write $N := H^p\gamma_3(H)$ (if $p = 3$ then $N = H^3$). Then $|H:N| = p^3$. On the other hand, since $G/K'$ is a group of maximal class of order at least $p^{p+2}$, its maximal subgroup $H/K'$ is also of maximal class and of order at least $p^{p+1}$, and this implies that $|H:K'N| = p^3$. Therefore $N = K'N$, i.e. $K' \le N$. Similarly, $L' \le N$, which leads to the contradiction $K'L' \le N < \gamma_3(G)$.

Finally, if $\gamma_3(G) \neq 1$, let $H = G/\gamma_4(G)$. Then $Z(H) = \gamma_3(G)/\gamma_4(G)$, and $H/Z(H)$ is one
of the two non-abelian groups of order $p^3$. But of these two, the one of exponent $p^2$ cannot be the central
factor group of any group (is {\it incapable}), therefore $H/Z(H)$ has exponent $p$, which is our claim.

\smallskip

For $p =3$ we have $p - 1 = 2$, hence

\proclaim {Proposition 5} Let $G$ be a 3-group such that $|G:G'| = 9$ and $|G| \ge 81$. Then $G^3 = \gamma_3(G)$ and that subgroup has index 27. All maximal subgroups of $G$ have at most three generators, and either one of them is metacyclic, or $G$ contains a maximal subgroup $H$ such that $|H:H'| \le 27$. If $|G:\gamma_4(G)| = 3^4$, then $G$ is of maximal class.

\endproclaim

{\bf Proof.} Most of this is either stated in Theorem 4, or was derived during the proof of part (f) there. The last claim is the case $p = 3$ of the fact, that if in a $p$-group $G$ the factor group $G/\gamma_{p+1}(G)$ is of maximal class, then so is $G$ [B1, 3.9].

\smallskip

\proclaim {Proposition 6} Let $G$ be a $p$-group, of order $p^5$ at least, in which $|G:G'| \le p^3$. Then all
maximal subgroups of $G$, with at most one exception, have irreducible characters of degree $p$. If an exceptional
maximal subgroup exists, then all other maximal subgroups $M$ satisfy $|M:M'| \le |G:G'|$.\endproclaim

{\bf Proof.} Since $|G:G'| \le p^3$, $G$ has irreducible characters of degree $p$. By [Ma1, Proposition 14(2)], if
some maximal subgroup $H$ does not have such characters, and $M$ is another maximal subgroup, then $|G':M'| \le
p$. That means that $|M:M'| \le |G:G'|$, and this shows that $M$ has irreducible characters of degree $p$. {\bf
QED}

\smallskip

An exceptional maximal subgroup may or may not exist. If $G$ is a $p$-group of maximal class, which is metabelian,
but in which the subgroup $G_1$ is not abelian, then all irreducible characters of all maximal
subgroups of $G$ have degrees 1 or $p$. On the other hand, in the groups of maximal class and order $p^{2r}$ ($p >
3,~2r = 6,8,...,p+1$) constructed by N.Blackburn [B1, pp. 61-62, or Hu, III.14.24], or the ones of order $p^{p+2}$
constructed by Slattery in [Sl2] ($p > 3$), it is not difficult to show that the maximal subgroups $H$ (in [B1])
or $E$ (in [Sl2]) have only irreducibles of degrees $1$ and $p^{r-1}$, or $1$ and $p^{\frac{p-1}{2}}$.

\smallskip

In the proofs of Theorems 1 to 3 we separate between groups of maximal class and others.

\smallskip

{\bf Proof of Theorem 1 for groups of maximal class, and of Theorem 2.} We may assume that $|G| \ge p^{p+2}$. By a famous result of Isaacs-Passman [Is1, 12.11], a $p$-group
$G$ has characters of degrees 1 and $p$ only, iff either $|G:Z(G)| \le p^3$ or $G$ contains an abelian maximal
subgroup. For groups of maximal class the first possibility means that $|G| \le p^4$, and then anyway $G$ has an
abelian maximal subgroup. Moreover, if $G$ has an abelian maximal subgroup, that subgroup must be $G_1$. Thus the
assumption that $G$ has characters of degree $deg(\chi) > p$ means that $G_1' \neq 1$. Now $[G_1,G_p] =
[G_1,G_1^p] = [G_1,G_1]^p$ is properly contained in $G_1'$. In $H = G/[G_1,G_p]$ we have $H_1' \neq 1$, therefore
$H$ has irreducibles of degree exceeding $p$. But $H_p \le Z(H_1)$, and $|H_1:H_p| = p^{p-1}$, therefore the
irreducible characters of $H_1$ have degrees at most $p^{\frac{p-1}{2}}$, and those of $H$ have degrees at most
$p^{\frac{p+1}{2}}$.

Next, let $G$ have characters of degree $deg(\chi) > p^2$. If $|G| \le p^8$, then $|G:Z(G)| \le p^7$, and all
irreducibles of $G$ have degree $p^3$ at most. Thus to prove our claim we may assume that $|G| \ge p^9$. According
to [Pa], a $p$-group $X$ has only characters of degree at most $p^2$, for $p$ odd, iff one of the following four
possibilities occur:

(i) $G$ contains an abelian subgroup of index $\le p^2$,

(ii) $G$ contains a maximal subgroup $H$ such that $|H:Z(H)| \le p^3$,

(iii) $|G:Z(G)| \le p^5$,

(iv) $|G:Z(G)| = p^6$, and if $H \ge Z(G)$ is a maximal subgroup of $G$, then $Z(H) = Z(G)$.

In a group of maximal class (iv) is impossible, because $C_G(Z_2(G))$ is a maximal subgroup, and (iii) means that
$|G| \le p^6$. (ii) means that either $G$ is an exceptional group of order at most $p^6$, or that $[G_1,G_4] = 1$,
and (i) means that $G'$ is abelian. Thus if $G$ has an irreducible character of degree at least $p^3$, then $G''
\neq 1 \neq [G_1,G_4]$. Note that by our assumption $G_8 \neq 1$, and then $G/G_7$ is non-exceptional, which
implies that $[G_1,G_4] \le G_6$. There are two indices $i,j,$ such that $G'' = G_i$ and $[G_1,G_4] = G_j$. First
assume that $i > j$, and let $H = G/G_{i+1}$. Then $H'' \neq 1 \neq [H_1,H_4]$. The last inequality shows that
$H_6 \neq 1$, therefore $|H| \ge p^7$. Thus $H$ violates (i)-(iv), and it has an irreducible character of degree
at least $p^3$. But $1 = (H'')^p = [H',(H')^p] = [H_2,H_{p+1}]$. Thus $H_{p+1} \le Z(H')$, and again the
characters of $H'$ have degrees at most $p^{\frac{p-1}{2}}$ and the characters of $H$ have degrees at most
$p^{\frac{p+3}{2}}$. If $i \le j$, we take $H = G/G_{j+1}$, and obtain $1 = [H_1,H_4]^p = [H_1,H_{p+3}]$, with
$|H_1:H_{p+3}| = p^{p+2}$, and proceed as before.

\smallskip

{\bf Proof of the rest of Theorem 1, and of Theorem 3.} Since $G$ is not of maximal class, neither is $G/\gamma_{p+1}(G)$ [Bl, 3.9]. Let $i \ge 2$ be the
first index such that $|\gamma_i(G)/\gamma_{i+1}(G)| > p$. Thus $3 \le i \le p$. Then
$|\gamma_i(G)/\gamma_{i+1}(G)| = p^2$. Theorem 4(a) shows that if $H$ is a maximal subgroup of $G$, then $H' \not\le
\gamma_{i+1}(G)$. The fact that $G$ has characters of degree $> p$ implies that $|G:Z(G)| > p^3$, therefore
$\gamma_3(G) \not\le Z(G)$, i.e. $\gamma_4(G) \neq 1$. Assume first that $i+1 \neq 4$, and let $K = G/\gamma_{i+1}(G)$. Then in $K = G/\gamma_{i+1}(G)$ no maximal subgroup is abelian, and the centre has index $> p^3$, therefore
$K$ has an irreducible character $\chi$ such that $deg(\chi)
> p$. On the other hand, $|K| = p^{i+2} \le p^{p+2}$, therefore $deg(\chi) \le p^{\frac{p+1}{2}}$.

If $i = 3$, then $\gamma_{i+1}(G) = \gamma_4(G) \neq 1$. Let $N \triangleleft G$ have index $p$ in
$\gamma_{i+1}(G)$, let $K = G/N$, and proceed as before.

If $|G:\gamma_4(G)| = p^5$, then we saw in Theorem 4 that if $C$ and $D$ are two maximal subgroups, then $C'
\not\le D'$. This implies that $C' \not\le G''$, and thus $L := G/G''$ does not have abelian maximal subgroups.
And $|L:Z(L)| > p^3$, because $[\gamma_3(G),G] = \gamma_4(G) \not\le G''$, and $\gamma_3(G)$ is the only normal
subgroup of index $p^3$. It follows that $G/G''$ has irreducible characters of degree bigger that $p$, and that
degree must be $p^2$, because $G'/G''$ is an abelian subgroup of index $p^2$. {\bf QED}

\bigskip

\centerline {\bf References}

\medskip
Be. Y.Berkovich, {\it Groups of Prime Power Order, vol. 1}, de Gruyter, Berlin 2008.

Bl. N.Blackburn, On a special class of $p$-groups, Acta Math. 100 (1958), 45-92.

Hu. B.Huppert, {\it Endliche Gruppen I}, Springer, Berlin 1967.

Is1. I.M.Isaacs, {\it Character Theory of Finite Groups}, Academic Press, San Diego 1976.

Is2. I.M.Isaacs, Sets of $p$-powers as irreducible character degrees, Proc. Amer. Math. Soc. 96 (1986),
551-552.

LGM. C.R.Leedham-Green and S.McKay, {\it The Structure of Groups of Prime Power Order}, Oxford University Press,
Oxford 2002.

Ma1. A.Mann, Minimal characters of $p$-groups, J. Gp. Th. 2 (1999), 225-250.

Ma2. A.Mann, More on normally monomial $p$-groups, in preparation.

Pa. D.S.Passman, Groups whose irreducible representations have degrees dividing $p^2$, Pac. J. Math. 17 (1966),
475-496.

Sl1. M.C.Slattery, Character degrees of normally monomial maximal class 5-groups, in {\it Character Theory of
Finite Groups (the Isaacs conference), Contemporary Mathematics 524}, American Mathematical Society, Providence
2010.

Sl2. M.C.Slattery, Maximal class $p$-groups with large character degree gaps, preprint.

\end